\documentclass[12pt]{article}
\usepackage{amssymb, amsmath, url, graphicx, hyperref}
\usepackage{indentfirst}

\def\3{\subset }
\def\4{\subseteq }
\def\<{\left<}
\def\>{\right>}

\def\bit{\begin{itemize}}
\def\eit{\end{itemize}}
\def\3{\subset }
\def\4{\subseteq }

\def\0{\leqno}

\def\barr{\begin{array}}
\def\earr{\end{array}}

\def\Z{{\rlap{$\kern2pt{\rm Z}$}{\rm Z}\,}}

\title{\bf On some probabilistic aspects of (generalized) dicyclic groups}
\author{Marius T\u arn\u auceanu and Mihai-Silviu Lazorec}
\date{December 6, 2016}

\begin{document}

\maketitle

\begin{abstract}
In this paper we study probabilistic aspects such as subgroup commutativity 
degree and cyclic subgroup commutativity degree of the (generalized) dicyclic groups.
We find explicit formulas for these concepts and we provide another example of a 
class of groups whose (cyclic) subgroup commutativity degree vanishes asymptotically.
\end{abstract}



\section{Introduction}

The link between finite groups and probability theories is a topic which is frequently 
studied in the last years. Given a finite group $G$, many results regarding the appartenence
of $G$ to a class of groups were obtained using a relevant concept known as the 
\textit{commutativity degree} of $G$ which measures the probability that two 
elements of the group commute. We refer the reader to \cite{2, 3} and \cite{5}-\cite{9} 
for more details. Inspired by this notion, in \cite{15,21}, we introduced another two probabilistic 
aspects called the \textit{subgroup commutativity degree} and the \textit{cyclic subgroup 
commutativity degree} of a finite group $G$. Denoting by $L(G)$ and $L_1(G)$ the subgroup 
lattice and the poset of cyclic subgoups of $G$, respectively, the above concepts 
are defined by 

$$sd(G)=\frac{1}{|L(G)|^2}\lbrace (H,K)\in L(G)^2 | HK=KH \rbrace=$$\\
$$ \hspace{0.9cm}=\frac{1}{|L(G)|^2}\lbrace (H,K)\in L(G)^2 | HK \in L(G)\rbrace$$ 

and

$$csd(G)=\frac{1}{|L_1(G)|^2}\lbrace (H,K)\in L_1(G)^2 | HK=KH \rbrace=$$\\
$$\hspace{1.3cm}=\frac{1}{|L_1(G)|^2}\lbrace (H,K)\in L_1(G)^2 | HK \in L_1(G)\rbrace,$$ respectively.

In other words, the (cyclic) subgroup commutativity degree measures the probability that 
two (cyclic) subgroups of $G$ commute. Moreover, for a subgroup $H$ of $G$, by introducing the sets 
$$C(H)=\lbrace K\in L(G)| HK=KH\rbrace$$ and $$C_1(H)=\lbrace K\in L_1(G)|HK=KH\rbrace,$$ we obtain 
$$sd(G)=\frac{1}{|L(G)|^2}\sum_{H\in L(G)}|C(H)|$$ and $$csd(G)=\frac{1}{|L_1(G)|}\sum_{H\in L_1(G)}|C_1(H)|.$$

Also, in \cite{15}-\cite{18} and \cite{21}, we 
provided explicit formulas for $sd(G)$ and $csd(G)$, when $G$ belongs to some particular classes of groups. Another purpose was to study the asymptotic behaviour of this notions. In this paper we will analyse the same problems regarding the dicyclic groups and some generalized dicyclic groups. 

For a positive integer $n \geqslant 1$, the dicyclic group, denoted by $Dic_{4n}$, is defined as 
$$Dic_{4n}= \langle a,\gamma | a^{2n}=e, \gamma^2=a^n, a^{\gamma}=a^{-1}\rangle.$$ This type of group 
has the following generalization: for an arbitrary abelian group $A$ of order $2n$, the generalized dicyclic group, denoted by $Dic_{4n}(A)$, is 
$$Dic_{4n}(A)=\langle A, \gamma | \gamma^4=e, \gamma^2\in A\setminus \lbrace e\rbrace, g^{\gamma}=g^{-1}, \forall g\in A\rangle.$$ 
\\

{\textbf{Remark 1.1.}} It is easy to see that $Dic_{4n}(\mathbb{Z}_{2n})=Dic_{4n}$. Regarding 
the study of the subgroup commutativity degree and of the cyclic subgroup commutativity degree 
of $Dic_{4n}(A)$, we mention that we will limit to the case $A\cong \mathbb{Z}_2 \times \mathbb{Z}_n$, since there is no explicit result which may allow us to compute the number of subgroups of an arbitrary abelian group of order $2n$.
\\

This paper is organised as follows. Section 2 provides explicit formulas for subgroup commutativity 
degree of the dicyclic group $Dic_{4n}$ and of the generalized dicyclic group $Dic_{4n}(\mathbb{Z}_2 \times \mathbb{Z}_n)$. For the last one, we cover a particular case, namely $n=2^m$, where $m\geqslant 2$ is a positive integer, which raise the problem of counting the number of subgroups of $\mathbb{Z}_2 \times D_{2^m}$ and of $\mathbb{Z}_2\times Q_{2^{m+1}}$. Results concerning the cyclic subgroup commutativity degree of the above mentioned groups are presented in Section 3. Some open problems and further research directions are indicated in the last section.   

Most of our notation is standard and will usually not be repeated here. Elementary notions and 
results on groups can be found in \cite{4,13}. For subgroup lattice concepts we refer the reader to
\cite{12,14,19}.           

\section{Subgroup commutativity degree of (generalized) dicyclic groups}

Our first purpose in this section is to compute explicitly the subgroup commutativity degree of the dicyclic group $Dic_{4n}$. Recall that, for a positive integer $n\geqslant 1$, the structure of this group is $$Dic_{4n}=\langle a,\gamma | a^{2n}=e, \gamma^2=a^n, a^{\gamma}=a^{-1}\rangle.$$   

To compute the subgroup commutativity degree of $Dic_{4n}$, it is helpful to know which are its subgroups. For each divisor $r$ of $2n$, the dicyclic group has one subgroup isomorphic to $\mathbb{Z}_r$, namely $H_0^r=\langle a^{\frac{2n}{r}}\rangle$. Also, for each divisor $s$ of $n$, $Dic_{4n}$ possesses $\frac{n}{s}$ subgroups isomorphic to $Dic_{4s}$, namely $H_i^s=\langle a^{\frac{n}{s}}, a^{i-1}\gamma\rangle$, where $i=\overline{1,\frac{n}{s}}$. 
\\

\textbf{Theorem 2.1.}
\textit{Let $n=2^mm'$, where $m\in \mathbb{N}$ and $m'$ is a positive odd number. Then the subgroup commutativity degree of the dicyclic group is 
$$sd(Dic_{4n})=\frac{(m+2)^2\tau(m')^2+2(m+2)\tau(m')\sigma(n)+[(m-1)2^{m+3}+9]g(m')}{[(m+2)\tau(m')+\sigma(n)]^2},$$
where $g$ is the arithmetic function given by 
$$g(m')=\prod\limits_{i=1}^k\frac{(2\alpha_i+1)p_i^{\alpha_i+2}-(2\alpha_i+3)p_i^{\alpha_i+1}+p_i+1}{(p_i-1)^2}$$ 
for the prime factorization of $m'$, i.e. $m'=p_1^{\alpha_1}p_2^{\alpha_2}\ldots p_k^{\alpha_k}.$}
\\

\textbf{Proof.} 
The subgroup commutativity degree of $Dic_{4n}$ is given by 
\begin{equation}\label{r1}
sd(Dic_{4n})=\frac{1}{|L(Dic_{4n})|^2}\bigg(\sum_{r|2n}|C(H_0^r)|+\sum_{s|n}\sum\limits_{i=1} ^\frac{n}{s}|C(H_i^s)|\bigg).
\end{equation}
Analysing the above mentioned subgroup structure of the dicyclic group, we infer that 
$$|L(Dic_{4n})|=\tau(2n)+\sigma(n)=(m+2)\tau(m')+\sigma(n).$$   
Let $r$ be a divisor of $2n$. The nature of the subgroup $H_0^r$ and the relation $a^{\gamma}=a^{-1}$, leads us to $$|C(H_0^r)|=|L(Dic_{4n})|=(m+2)\tau(m')+\sigma(n).$$ 
Hence $$\sum_{r|2n}|C(H_0^r)|=(m+2)^2\tau(m')^2+(m+2)\tau(m')\sigma(n).$$ 
Now, let $s$ be a divisor of $n$ and consider the dicylic subgroup $H_i^s$ of $Dic_{4n}$, where $i\in \lbrace 1,2,\ldots, \frac{n}{s}\rbrace$. Then $$C(H_i^s)=\bigg( \bigcup_{r|2n}\lbrace H_0^r\rbrace \bigg) \cup \bigg\lbrace H_j^t | \ H_i^sH_j^t=H_j^tH_i^s, where \ t|n, j\in \lbrace 1, 2, \ldots, \frac{n}{t} \rbrace \bigg\rbrace.$$ 
It is known that $Z(Dic_{4n})=\langle a^n\rangle$ and $Z(Dic_{4n})\subset H_i^s$, for each divisor $s$ of $n$ and for all $i=\overline{1,\frac{n}{s}}$. Also, it is easy to observe that $\frac{Dic_{4n}}{Z(Dic_{4n})}\cong D_{2n},$ where $D_{2n}=\langle x,y | x^n=y^2=e, x^y=x^{-1}\rangle$ is the dihedral group of order $2n$. Thus 
$$H_j^t \in C(H_i^s) \Leftrightarrow H_i^sH_j^t=H_j^tH_i^s \Leftrightarrow \frac{H_i^s}{Z(Dic_{4n})}\frac{H_j^t}{Z(Dic_{4n})}=\frac{H_j^t}{Z(Dic_{4n})}\frac{H_i^s}{Z(Dic_{4n})}.$$
Notice that, by the isomorphism $\frac{Dic_{4n}}{Z(Dic_{4n})}\cong D_{2n}$, the subgroups $\frac{H_i^s}{Z(Dic_{4n})}$ and $\frac{H_j^t}{Z(Dic_{4n})}$ of $\frac{Dic_{4n}}{Z(Dic_{4n})}$ are isomorphic with the subgroups $K_i^s=\langle x^\frac{n}{s}, x^{i-1}y \rangle$ and $K_j^t=\langle x^\frac{n}{t}, x^{j-1}y\rangle$ of $D_{2n}$. Hence, the above equivalence actually says that the subgroups $H_i^s$ and $H_j^t$ of $Dic_{4n}$ commute iff the subgroups $K_i^s$ and $K_j^t$ of $D_{2n}$ commute. This problem was studied in \cite{15} and, according to this reference, we obtain
$$\sum_{s|n}\sum\limits_{i=1} ^\frac{n}{s}|C(H_i^s)|=(m+2)\tau(m')\sigma(n)+[(m-1)2^{m+3}+9]g(m').$$  

Since we computed all the unknown quantities that appear in \eqref{r1}, our proof is complete. 
\hfill\rule{1,5mm}{1,5mm}
\\

\textbf{Example 2.2.}

$\bullet sd(Dic_4)=sd(\mathbb{Z}_4)=1;$\\

$\bullet sd(Dic_8)=sd(Q_8)=1;$\\

$\bullet sd(Dic_{12})=\frac{29}{32};$\\

$\bullet sd(Dic_{16})=\frac{113}{121};$
\\

We continue our study by finding the subgroup commutativity degree of the generalized dicyclic group $Dic_{4n}(A)$, where $A=\mathbb{Z}_2\times \mathbb{Z}_n$ and $n$ is a positive even integer. We recall that the structure of this group is 
$$Dic_{4n}(A)=\langle A, \gamma | \gamma^4=e, \gamma^2\in A\setminus \lbrace e\rbrace, g^{\gamma}=g^{-1}, \forall g\in A\rangle.$$

Let $a$ and $b$ be the generators of the cyclic groups $\mathbb{Z}_n$ and $\mathbb{Z}_2$, respectively. Since $\gamma^2\in A\setminus \lbrace e\rbrace$ and $\gamma^4=e$, we infer that $\gamma^2\in \lbrace a^{\frac{n}{2}}, b, a^\frac{n}{2}b \rbrace$. 
In the proof of the main result which provides the explicit formula of $sd(Dic_{4n}(A))$ we will need some preliminary results to analyse the case characterised by $n=2^m$, where $m\geqslant 2$ is a positive integer and $\gamma^2=a^{\frac{n}{2}}$. Denote by $n'$ the quantity $n'=\frac{n}{2}=2^{m-1}$. The necessary tools are given by the following lemma and corollary.
\\

\textbf{Lemma 2.3.} \textit{The number of subgroups of the direct product $\mathbb{Z}_2 \times D_{2n'}$ is given by
$$|L(\mathbb{Z}_2 \times D_{2n'})|=5\sigma(n')+3\tau(n')-2n'-1.$$}

\textbf{Proof.}
We will apply Goursat's Lemma to prove the desired equality. This result states that there is a bijection between the set of subgroups of $\mathbb{Z}_2 \times D_{2n'}$ and the set of $5$-uples $(H_1,H_1',H_2,H_2',\phi)$, where $H_1 \triangleleft H_1' \leqslant \mathbb{Z}_2$, $H_2 \triangleleft H_2' \leqslant D_{2n'}$ and $\phi:\frac{H_1'}{H_1}\longrightarrow \frac{H_2'}{H_2}$ is an isomorphism. Moreover, a subgroup $H$ of $\mathbb{Z}_2\times D_{2n'}$ is $H=\lbrace (g_1,g_2)\in \mathbb{Z}_2\times D_{2n'}| \phi(g_1H_1)=g_2H_2\rbrace$.

Since $|L(D_{2n'})|=\tau(n')+\sigma(n')$, there are $2(\tau(n')+\sigma(n'))$ subgroups of $\mathbb{Z}_2\times D_{2n'}$ corresponding to the $5$-uples $(H_1,H_1,H_2,H_2,\phi)$. For each divisor $d\not=n'$ of $n'$, we get a $5$-uple $(\lbrace 1\rbrace, \mathbb{Z}_2,\langle x^\frac{2n'}{d}\rangle, \langle x^\frac{n'}{d}\rangle,\phi)$, where $\lbrace 1\rbrace$ denotes the trivial subgroup of $\mathbb{Z}_2$. Hence, we obtain other $\tau(n')-1$ subgroups of $\mathbb{Z}_2 \times D_{2n'}$. Going further, it is easy to remark that $n$ subgroups of our direct product are in bijection with the $5$-uples $(\lbrace 1\rbrace, \mathbb{Z}_2, \lbrace e\rbrace, \langle x^{i-1}y\rangle, \phi)$, $i=\overline{1,n}$.

To count the number of the remaining subgroups of $\mathbb{Z}_2 \times D_{2n'}$, we take $d\not=1$ to be a divisor of $n'$. Remark that for obtaining the isomorphism $\phi$, a dihedral subgroup of $D_{2n'}$, namely $\langle x^\frac{n'}{d}, x^{i-1}y\rangle$, where $i\in \lbrace 1,2,\ldots \frac{n'}{d}\rbrace$ is fixed, could be factored only by $\langle x^{\frac{n'}{d}}\rangle$ and $\langle x^\frac{2n'}{d}, x^{j-1}y\rangle$, where $j\in \lbrace i, i+\frac{n'}{d}\rbrace$. Consequently we have 
$$\sum\limits_{\substack{d|n' \\ d\not=1}}\sum\limits_{i=1}^\frac{n'}{d}=3(\sigma(n')-n')$$
additional subgroups of our direct product. By inspecting the lattice of subgroups of the dihedral group $D_{2n'}$, we specify that, using its subgroups, we can not have other posibilities to form a factor group whose order is less or equal with two, meaning that we can not find other $5$-uples. Our proof is complete once we add up the above obtained values. 
\hfill\rule{1,5mm}{1,5mm}
\\

Recall that the generalized quaternion group has the following structure 
$$Q_{4n'}=\langle x,y | x^n=y^4=1, x^y=x^{-1}\rangle.$$
It is known that $Q_{4n'}$ has an unique minimal subgroup, which is its center $Z(Q_{4n'})=\langle x^{n'}\rangle$, and $\frac{Q_{4n'}}{Z(Q_{4n'})}\cong D_{2n'}$. We infer that 
$$\frac{\mathbb{Z}_2 \times Q_{4n'}}{Z(Q_{4n'})}\cong \mathbb{Z}_2 \times D_{2n'}.$$ 
Using the Correspondence Theorem and the fact that there are $3$ subgroups of $\mathbb{Z}_2 \times Q_{4n'}$ which do not contain $Z(Q_{4n'})$, namely the trivial subgroup, $\mathbb{Z}_2$ and $\langle x^{n'}b\rangle$, we obtain the following result.
\\

\textbf{Corollary 2.4.}
\textit{The number of subgroups of the group $\mathbb{Z}_2 \times Q_{4n'}$ is 
given by 
$$|L(\mathbb{Z}_2\times Q_{4n'})|=5\sigma(n')+3\tau(n')-2n'+2.$$}
\\

\textbf{Remark 2.5.} Since $n'=2^{m-1}$, the number of divisor of $n'$ and their sum are $\tau(n')=m$ and $\sigma(n')=2^m-1$, respectively. Then the conclusion of \textbf{Corollary 2.4.} can be written as
$$|L(\mathbb{Z}_2\times Q_{2^{m+1}})|=2^{m+2}+3(m-1).$$ 
\\
       
We are ready to prove the main result concerning the subgroup commutativity degree of the generalized dicyclic group $Dic_{4n}(\mathbb{Z}_2 \times \mathbb{Z}_n)$.
\\

\textbf{Theorem 2.6.} \textit{Let $A=\mathbb{Z}_2\times \mathbb{Z}_n$ be an abelian group, where $n=2^mm'$, $m\in\mathbb{N}^*$ and $m'$ is a positive odd integer. Then, \\ \\
(i) \  if $m\geqslant 2, m'\not=1$ and $\gamma^2\in\lbrace a^{\frac{n}{2}}, b, a^{\frac{n}{2}}b\rbrace$ or if $m\geqslant 2, m'=1$ and $\gamma^2 \in \lbrace b, a^\frac{n}{2}b \rbrace$, the subgroup commutativity degree of the generalized dicyclic group $Dic_{4n}(A)$ is 
$$sd(Dic_{4n}(A))=\frac{|L(A)|^2+2|L(A)|\sigma (n)+[(m-1)2^{m+3}+9]g(m')}{(|L(A)|+\sigma (n))^2},$$   where $g$ is the arithmetic function given by 
$$g(m')=\prod\limits_{i=1}^k\frac{(2\alpha_i+1)p_i^{\alpha_i+2}-(2\alpha_i+3)p_i^{\alpha_i+1}+p_i+1}{(p_i-1)^2}$$ 
for the prime factorization of $m'$, i.e. $m'=p_1^{\alpha_1}p_2^{\alpha_2}\ldots p_k^{\alpha_k}$.\\
(ii) \ if $m\geqslant 2, m'=1$ and $\gamma^2=a^{\frac{n}{2}}$, the generalized dicyclic group $Dic_{4n}(A)$ is isomorphic to the direct product $\mathbb{Z}_2\times Q_{2^{m+1}}$ and its subgroup commutativity degree is
$$sd(\mathbb{Z}_2\times Q_{2^{m+1}})=\frac{2^{m+2}(24m-37)+9m^2-18m+185}{[2^{m+2}+3(m-1)]^2}.$$\\
(iii) if $m=m'=1$ and $\gamma^2\in \lbrace a^\frac{n}{2}, b, a^\frac{n}{2}b\rbrace$, the generalized dicyclic group $Dic_{4n}(A)$ is isomorphic to the abelian group $\mathbb{Z}_2^3$ and its subgroup commutativity degree is $$sd(\mathbb{Z}_2^3)=1.$$ }
\\

\textbf{Proof.} (i) Let $n=2^mm'$, where $m\geqslant 2$, $m'$ is a positive odd integer and $\gamma^2 \in \lbrace b, a^\frac{n}{2}b\rbrace$. In this case, besides the subgroups of the abelian group $A$, for each divisor $r$ of $n$, the generalized dicyclic group $Dic_{4n}(A)$ contains the subgroups $H_i^r=\langle a^\frac{n}{r}, a^{i-1}\gamma\rangle$, where $i=\overline{1,\frac{n}{r}}$. Hence 
$$|L(Dic_{4n}(A))|=|L(A)|+\sigma(n).$$

Since $g^\gamma=g^{-1}, \forall g\in A$ and $A$ is an abelian group, it is easy to see that the subgroups of $A$ commute with all subgroups of $Dic_{4n}(A)$. Consequently there are $|L(A)|(|L(A)|+\sigma(n))$ pairs of subgroups of the generalized dicyclic group which commute. Hence 
\begin{equation}\label{r3}
sd(Dic_{4n}(A))=\frac{|L(A)|(|L(A)|+\sigma(n))+\sum\limits_{s|n}\sum\limits_{i=1} ^\frac{n}{s}|C(H_i^s)|}{(|L(A)|+\sigma(n))^2}.
\end{equation}
Moreover, for a divisor $r$ of $n$ and a fixed $i\in \lbrace 1, 2, \ldots, \frac{n}{r}\rbrace$, we have
$$C(H_i^r)=L(A)\cup \lbrace H_j^t | H_i^rH_j^t=H_j^tH_i^r, \ t|n \ and \ j\in \lbrace 1, 2, \ldots, \frac{n}{t}\rbrace\rbrace.$$

Denote by $x, y$ and $e_1$ the elements $a\langle \gamma^2\rangle, \gamma \langle\gamma^2\rangle$ and $e\langle\gamma^2\rangle$, respectively. Then 
$$\frac{Dic_{4n}(A)}{\langle \gamma^2\rangle}=\langle x,y| x^n=y^2=e_1, x^y=x^{-1}\rangle\cong D_{2n}.$$  
Using the same reasoning as in the proof of \textbf{Theorem 2.1.}, the subgroups $H_i^r$ and $H_j^t$ of $Dic_{4n}(A)$ commute iff the corresponding subgroups, through the above isomorphism, of $D_{2n}$ commute. We obtain
$$\sum_{s|n}\sum\limits_{i=1} ^\frac{n}{s}|C(H_i^s)|=|L(A)|\sigma(n)+[(m-1)2^{m+3}+9]g(m').$$
Using the last relation in \eqref{r3}, we finish the proof in this case.

If $m\geqslant 2$, $m'\not=1$ is a positive odd number and $\gamma^2=a^\frac{n}{2}$, besides the subgroups of $A$, for each divisor $r$ of $n$, the generalized dicyclic group $Dic_{4n}(A)$ contains the subgroups $\langle(ab)^\frac{n}{r}, (ab)^{i-1}\gamma\rangle$, where $i=\overline{1,\frac{n}{r}}$. The steps of the proof are similar with the ones we already made and we get the same explicit formula for the subgroup commutativity degree of $Dic_{4n}$. Hence, the statement (i) is proven.    
\\

(ii) Let $m\geqslant 2, m'=1$ and $\gamma^2=a^\frac{n}{2}$. The generalized dicyclic group in this case is 
$$Dic_{4n}(A)=\langle a, b, \gamma | a^n=b^2=\gamma^4=e, ab=ba, a^\gamma=a^{-1}, b^\gamma=b^{-1}\rangle.$$
Again, the subgroups of the abelian group $A$ commute with all subgroups of $Dic_{4n}(A)$ as a consequence of the equality $g^\gamma=g^{-1}, \forall g\in A$. Also, by the same reason, the subgroups of $Dic_{4n}(A)$ which are not contained in $A$ commute with the subgroups of $A$.
Since $\langle b\rangle \cap \langle a,\gamma |a^n=\gamma^4=e, a^\gamma=a^{-1}\rangle=\lbrace e\rbrace$ and the same subgroups commute because $ab=ba$ and $b^\gamma=b^{-1}$, we infer that 
$$Dic_{4n}(A)\cong \mathbb{Z}_2\times Q_{2^{m+1}}.$$
Hence, finding an explicit formula for the generalized dicyclic group $Dic_{4n}(A)$ is equivalent to the same problem related to the direct product $\mathbb{Z}_2\times Q_{2^{m+1}}$. Until now, we proved that there are $|L(A)|(2|L(\mathbb{Z}_2\times Q_{2^{m+1}})|-|L(A)|)$ pairs of subgroups which commute. Further,  we need to count the pairs of subgroups of $(L(\mathbb{Z}_2\times Q_{2^{m+1}})\setminus L(A))^2$ with the same property.

Denote by $x, y, b_1$ and $e_1$ the elements $a\langle \gamma^2\rangle, \gamma\langle \gamma^2\rangle, b\langle \gamma^2\rangle$ and $e\langle \gamma^2\rangle$, respectively. We get the following isomorphism 
$$\frac{\mathbb{Z}_2\times Q_{2^{m+1}}}{\langle \gamma^2\rangle}=\langle b_1\rangle\times\langle x, y| x^{2^m}=y^2=e_1, x^y=x^{-1}\rangle\cong \mathbb{Z}_2\times D_{2n'},$$ where $n'=\frac{n}{2}=2^{m-1}.$
All subgroups of $L(\mathbb{Z}_2\times Q_{2^{m+1}})\setminus L(A)$ contain $\langle \gamma^2\rangle$, so their commutativity is equivalent with the commutativity of the corresponding subgroups of the direct product $\mathbb{Z}_2\times D_{2n'}$ through the above isomorphism. 

These subgroups are divided into the following 5 sets:
$$H_1=\bigg\lbrace\langle x^\frac{n'}{r}, x^{i-1}y\rangle| \ r|n', i=\overline{1,\frac{n'}{r}}\bigg\rbrace,$$
$$H_2=\bigg\lbrace\langle x^\frac{n'}{r}b_1, x^{i-1}y\rangle| \ r|n', i=\overline{1,\frac{n'}{r}}\bigg\rbrace,$$
$$H_3=\bigg\lbrace\langle x^\frac{n'}{r}, x^{i-1}yb_1\rangle| \ r|n', i=\overline{1,\frac{n'}{r}}\bigg\rbrace,$$
$$H_4=\bigg\lbrace \langle x^\frac{n'}{r}b_1, x^{i-1}yb_1\rangle| \ r\not=1, r|n', i=\overline{1,\frac{n'}{r}}\bigg\rbrace,$$
$$H_5=\bigg\lbrace \langle b_1,x^\frac{n'}{r}, x^{i-1}y\rangle| \ r\not=1, r|n', i=\overline{1,\frac{n'}{r}}\bigg\rbrace.$$
The general form of one of the above listed subgroups is 
$$H_i^r(\alpha,\beta,\delta)=\langle b_1^\alpha, x^\frac{n'}{r}b_1^\beta, x^{i-1}yb_1^\delta\rangle,$$ 
where $r|n', i=\overline{1,\frac{n'}{r}}$ and $\alpha,\beta, \delta\in \lbrace 0, 1\rbrace$. Another useful way to write $H_i^r$ is the following one 
$$H_i^r(\alpha,\beta,\delta)=H_0^r(\alpha,\beta)\langle x^{i-1}yb_1^\delta\rangle,$$
where $H_0^r(\alpha,\beta)=\langle b_1^\alpha, x^\frac{n'}{r}b_1^\beta\rangle$ is a normal subgroup of $\mathbb{Z}_2 \times D_{2n'}$ contained in $\langle b_1,x\rangle\cong \mathbb{Z}_2\times \mathbb{Z}_{n'}$. Taking another subgroup of the same form, namely $H_j^s(\alpha_1,\beta_1, \delta_1)=\langle b_1^{\alpha_1}, x^\frac{n'}{s}b_1^{\beta_1}, x^{j-1}yb_1^{\delta_1}\rangle$, where $s|n', j=\overline{1,\frac{n'}{s}}, \alpha_1, \beta_1, \delta_1 \in \lbrace 0,1 \rbrace$, and noticing that $b_1\in Z(\mathbb{Z}_2\times D_{2n'})$, we obtain
$$H_i^r(\alpha,\beta,\delta)H_j^s(\alpha_1,\beta_1, \delta_1)=H_j^s(\alpha_1,\beta_1, \delta_1)H_i^r(\alpha,\beta,\delta) \Leftrightarrow$$
\begin{equation}\label{r4}
\Leftrightarrow x^{2(i-j)}\in H_0^r(\alpha,\beta)H_0^s(\alpha_1,\beta_1).
\end{equation}  

The structure of $H_0^r(\alpha,\beta)H_0^s(\alpha_1,\beta_1)$ is changing when we choose two arbitrary subgroups from the sets $H_k, k=\overline{1,5}$, but we can determine it by inspecting the lattice of $\mathbb{Z}_2\times \mathbb{Z}_n'$. We denote by $c_{k_1k_2}$ the number of pairs of subgroups which commute contained in the sets $H_{k_1} \times H_{k_2}, k_1,k_2=\overline{1,5}$ with $k_2\geqslant k_1$ (since we have a simmetry). Our next purpose is to compute these quantities. We mention that $k_1,k_2\in\lbrace 4,5\rbrace$ implies $r\not=1$ and $s\not=1$. \\

$(c1) \ k_1=k_2=1, k_1=k_2=3$ and $k_1=1, k_2=3$

The relation \eqref{r4} becomes 
$$x^{2(i-j)}\in H_0^r(0,0)H_0^s(0,0)=\langle x^{\frac{n'}{[r,s]}} \rangle \Leftrightarrow \frac{n'}{[r,s]}|2(i-j),$$ where $i=\overline{1,\frac{n'}{r}}$ and $j=\overline{1,\frac{n'}{s}}.$
In \cite{15}, it was showed that the above relation has $2^{m+2}(m-2)+9$ solutions. Hence, 
$$c_{k_1k_2}=2^{m+2}(m-2)+9.$$  \\

$(c2) \ k_1=k_2=2$

In this case, \eqref{r4} is written as
$$x^{2(i-j)}\in H_0^r(0,1)H_0^s(0,1)=\begin{cases} \langle b_1, x^\frac{n'}{[r,s]}\rangle &\mbox{, \ if } r\not=s \\  \langle x^\frac{n'}{r}b_1\rangle &\mbox{, \ if } r=s \end{cases}.$$

If $r=s=1$, then $i\equiv j(mod \ \frac{n'}{2})$, where $i, j=\overline{1,n'}$. For each $i$, the congruence has 2 solutions $j\in \lbrace 1, 2, \ldots, n'\rbrace$. If $r=s>1$, then $i\equiv j(mod \ \frac{n'}{r})$, where $i, j=\overline{1,\frac{n'}{r}}$. For every $i$, the congruence has only one solution. We get a total of  
$$2n'+\sum\limits_{\substack{r|n' \\ r>1}}\frac{n'}{r}=n'+\sigma(n')$$ pairs of subgroups which commute corresponding to the case $r=s$.    

If $r\not=s$, we must count the number of solutions of the relation $\frac{n'}{[r,s]}|2(i-j).$
If $r=s=n'$, we find only one solution. If $r=s<n'$, the relation is equivalent to $i\equiv j( mod \  \frac{n'}{2r})$, where $i, j=\overline{1,\frac{n'}{r}}$. For each $i$, the previuous congruence has 2 solutions. Hence, for $r=s$, we have 
$$1+\sum\limits_{\substack{r|n' \\ r\not=n'}}2\frac{n'}{r}=2\sigma(n')-1$$ solutions. Consequently, to find the pairs of subgroups which commute in the case $r\not=s$, we must substract the above number from $2^{m+2}(m-2)+9$. We obtain 
$$c_{k_1k_2}=n'+\sigma(n')+2^{m+2}(m-2)+9-2\sigma(n')+1=2^{m+2}(m-2)-2^{m-1}+11.$$  
\\

$(c3) \ k_1=1, k_2=2$ and $k_1=2, k_2=3$

Now, \eqref{r4} becomes
$$x^{2(i-j)}\in H_0^r(0,0)H_0^s(0,1)=\begin{cases} \langle x^\frac{n'}{s}b_1\rangle &\mbox{, \ if } r<s \\ \langle b_1, x^\frac{n'}{[r,s]}\rangle &\mbox{, \ if } r\geqslant s  \end{cases}.$$  

Let $r=2^u$ and $s$ be divisors of $n'$ such that $r<s$. We infer that $u\in \lbrace 0, 1, \ldots, m-2\rbrace$ and  $i\equiv j (mod \ \frac{n'}{s})$, where $i=\overline{1,\frac{n'}{r}}$ and $j=\overline{1,\frac{n'}{s}}$. For each $i$ the congruence has one solution. There are $\tau(n')-(u+1)$ divisors $s$ of $n'$ which are greater than the divisor $r$ of $n'$. Consequently, we get 
$$\sum\limits_{u=0}^{m-2}{2^{m-1-u}(\tau(n')-(u+1))}=\sum\limits_{t=1}^{m-1}2^tt=2^m(m-2)+2$$
pairs of subgroups which commute corresponding to the case $r<s$. 

We already saw that for $r=s$, the relation $\frac{n'}{[r,s]}|2(i-j)$ has $2\sigma(n')-1$ solutions. Now, we need to discuss the case $r>s$. If $r=n'$, then we obtain 
$$\sum\limits_{\substack{s|n' \\ s\not=n'}}\sum\limits_{j=1}^\frac{n'}{s} 1=\sigma(n')-1$$
more solutions of the same relation. Let $r=2^u$ and $s=2^v$ be divisors of $n'$ such that $r>s$ and $r\not=n'$. It follows that $u\in \lbrace 1,2, \ldots, m-2\rbrace$ and $v\in \lbrace 0, 1, \ldots, m-3\rbrace$. We can rewrite our relation as $i\equiv j (mod \ 2^{m-2-u})$. If $r$ is fixed, for each $i$ this congruence has 
$$\sum\limits_{v=0}^{u-1}2^{u-v+1}=2^{u+2}\bigg[ 1-\bigg(\frac{1}{2}\bigg)^u\bigg]$$ 
solutions. Consequently, if $r\geqslant s$, we obtain
$$3\sigma(n')-2+\sum\limits_{u=1}^{m-2} 2^{m+1}\bigg[ 1-\bigg(\frac{1}{2}\bigg)^u\bigg]=3\sigma(n')-2+2^{m+1}[(m-2)+2^{2-m}-1]$$
pairs of subgroups which commute in the case $r\geqslant s$. We coclude that 
$$c_{k_1k_2}=3\sigma(n')+2^m(m-2)+2^{m+1}[(m-2)+2^{2-m}-1]=2^m(3m-5)+5.$$
\\

$(c4) \ k_1=1, k_2=4$ and $k_1=3, k_2=4$

Then \eqref{r4} states that
$$x^{2(i-j)}\in H_0^r(0,0)H_0^s(0,1)=\begin{cases} \langle x^\frac{n'}{s}b_1\rangle &\mbox{, \ if } r<s \\ \langle b_1, x^\frac{n'}{[r,s]}\rangle &\mbox{, \ if } r\geqslant s  \end{cases}.$$
The steps in this case are similar with the ones we made in $(c3)$, the only change being that $s\not=1$. We obtain
$$c_{k_1k_2}=3(\sigma(n')-n')+2^m(m-2)+2^m[(m-2)+2(2^{2-m}-1)]=2^{m-1}(4m-9)+5.$$ 
\\

$(c5) \ k_1=1, k_2=5$ and $k_1=3, k_2=5$

Relation \eqref{r4} can be written as 
$$x^{2(i-j)}\in H_0^r(0,0)H_0^s(1,0)=\langle b_1, x^{\frac{n'}{[r,s]}}\rangle,$$
which is equivalent to counting the solutions of the relation $\frac{n'}{[r,s]}|2(i-j)$, where $r$ and $s$ are divisors of $n', s\not=1$ and $i=\overline{1,\frac{n'}{r}}, j=\overline{1,\frac{n'}{s}}$. If $r=n'$, our relation has 
$$\sum\limits_{\substack{s|n' \\ s\not=1}}\sum\limits_{j=1}^{\frac{n'}{s}}1=\sigma(n')-n'$$ solutions.

If $s=n'$ and $r\not=n'$ is fixed the above relation has only one solution for each $i$. Take $r=2^u$ and $s=2^v$, where $u=\overline{0,m-2}$ and $v=\overline{1,m-2}$. Then, the divisibility is equivalent to the congruence $i\equiv j(mod \ 2^{max\lbrace u,v \rbrace-v+1})$. For a fixed $r$, this congruence has 
$$1+\sum\limits_{v=1}^{m-2}2^{max\lbrace u,v \rbrace-v+1}=2^{u+1}-2u+2m-5$$
solutions for each $i$. Hence, the total number of pairs of subroups which commute is
$$c_{k_1k_2}=\sigma(n')-n'+\sum\limits_{u=0}^{m-2}2^{m-1-u}(2^{u+1}-2u+2m-5)=\sigma(n')-n'+2^m(3m-8)+10=$$
$$\hspace{1cm} =2^{m-1}(6m-15)+9.$$       
\\ 

$(c6) \ k_1=2, k_2=4$

The relation \eqref{r4} becomes 
$$x^{2(i-j)}\in H_0^r(0,1)H_0^s(0,1)=\begin{cases} \langle b_1, x^\frac{n'}{s}\rangle &\mbox{, \ if } r=1\\ \langle b_1, x^\frac{n'}{[r,s]}\rangle &\mbox{, \ if } r\not=s \ and \ r\not=1 \\  \langle x^\frac{n'}{s}b_1\rangle &\mbox{, \ if } r=s \end{cases}.$$

We will use some results that we already proved and edit them since $s\not=1$. By $(c2)$, if $r=s$, we have $\sigma(n')-n'$ pairs of subgroups which commute. For $r\not=s$ and $r\not=1$, by inspecting $(c3)$ and $(c4)$, we get 
$$2\lbrace(\sigma(n')-n'-1)+2^m[(m-2)+2(2^{2-m}-1)]\rbrace$$
more solutions. Let $r=1$. Then, we obtain the relation $\frac{n'}{s}|2(i-j)$, where $i=\overline{1,n'}$ and $j=\overline{1,\frac{n'}{s}}$, which provides $n'+2n'(m-2)$ pairs of subgroups which commute. Adding up the obtained number, it follows that
$$c_{k_1k_2}=3\sigma(n')+2n'(m-3)-2+2^{m+1}[(m-2)+2(2^{2-m}-1)]=2^m(3m-8)+11.$$
\\

$(c7) \ k_1=2, k_2=5$

In this case, \eqref{r4} states that
$$x^{2(i-j)}\in H_0^r(0,1)H_0^s(1,0)=\begin{cases} \langle b_1, x^\frac{n'}{s}\rangle &\mbox{, \ if } r=1\\ \langle b_1, x^\frac{n'}{[r,s]}\rangle &\mbox{, \ if } r\not=1 \end{cases}.$$  

Again, we will use our work done in the previous cases. Let $r\not=1$. By $(c2)$, if $r=s$, the relation $\frac{n'}{[r,s]}|2(i-j)$ has $1+2(\sigma(n')-n'-1)$. Using $(c3)$ and $(c4)$, for $r\not=s$, the same relation has 
$$2\lbrace (\sigma(n')-n'-1)+2^m[(m-2)+2(2^{2-m}-1)]\rbrace$$ more solutions. Also, $(c6)$ provides the number of pairs of subgroups which commute if $r=1$, this being $n'+2n'(m-2).$ Hence 
$$c_{k_1k_2}=1+4(\sigma(n')-n'-1)+n'+2n'(m-2)+2^{m+1}[(m-2)+2(2^{2-m}-1)]=$$
$$\hspace{1cm} =2^{m-1}(6m-15)+9.$$
\\

$(c8) \ k_1=k_2=4$

Relation \eqref{r4} becomes
$$x^{2(i-j)}\in H_0^r(0,1)H_0^s(0,1)=\begin{cases} \langle x^\frac{n'}{r}b_1\rangle &\mbox{, \ if } r=s\\ \langle b_1, x^\frac{n'}{[r,s]}\rangle &\mbox{, \ if } r\not=s \end{cases}.$$         

By $(c2), (c3)$ and $(c4)$, we obtain
$$c_{k_1k_2}=\sigma(n')-n'+2(\sigma(n')-n'-1)+2^{m+1}[(m-2)+2(2^{2-m}-1)]=$$
$$\hspace{1.4cm} =2^{m-1}(4m-13)+11.$$
\\

$(c9) \ k_1=4, k_2=5$ and $k_1=k_2=5$

In this case, \eqref{r4} becomes
$$x^{2(i-j)}\in H_0^r(0,1)H_0^s(1,0)=\langle b_1, x^\frac{n'}{[r,s]}\rangle.$$

Using the same arguments as in $(c7)$, the number of pairs of subgroups which commute is 
$$c_{k_1k_2}=1+4(\sigma(n')-n'-1)+2^{m+1}[(m-2)+2(2^{2-m}-1)]=2^m(2m-6)+9.$$ 
\\

Consequently, the number of pairs of subgroups which commute contained in the sets $H_{k_1}\times H_{k_2}$, where $k_1, k_2=\overline{1,5}$, is 
$$4c_{11}+c_{22}+4c_{12}+4c_{14}+4c_{15}+2c_{24}+2c_{25}+c_{44}+3c_{45}=2^m(72m-164)+201.$$ 
Hence, the subgroup commutativity degree of $\mathbb{Z}_2 \times Q_{2^{m+1}}$ is given by 
$$sd(\mathbb{Z}_2\times Q_{2^{m+1}})=\frac{|L(A)|(2|L(\mathbb{Z}_2\times Q_{2^{m+1}})|-|L(A)|)+2^m(72m-164)+201}{|L(\mathbb{Z}_2\times Q_{2^{m+1}})|^2}.$$
Recall that the subgroups forming the sets $H_k, k=\overline{1,5}$ are subgroups of $\mathbb{Z}_2\times Q_{4n'}$ which are not contained in the abelian group $A$. Their number is $5\sigma(n')-2n'$. \textbf{Corollary 2.4} leads us to $|L(A)|=3\tau(n')+2=3m+2$. Since the total number of subgroups of 
$\mathbb{Z}_2\times Q_{2^{m+1}}$ is also known by \textbf{Remark 2.5}, we obtain the explicit formula
$$sd(\mathbb{Z}_2\times Q_{2^{m+1}})=\frac{2^{m+2}(24m-37)+9m^2-18m+185}{[2^{m+2}+3(m-1)]^2}.$$
\\

(iii) Let $m=m'=1$ and $\gamma^2\in \lbrace a^\frac{n}{2}, b, a^\frac{n}{2}b\rbrace$. Since, in this case, $Dic_{4n}(A)$ is an abelian group of order $8$ generated by $3$ elements of order $2$, the conclusion follows easily. 
\hfill\rule{1,5mm}{1,5mm}
\\

\textbf{Example 2.7.}

$\bullet \ sd(Dic_{4n}(\mathbb{Z}_2\times \mathbb{Z}_4))=sd(\mathbb{Z}_2\times Q_8)=1;$\\

$\bullet \ sd(Dic_{4n}(\mathbb{Z}_2\times \mathbb{Z}_6))=\frac{215}{242};$\\

$\bullet \ sd(Dic_{4n}(\mathbb{Z}_2\times \mathbb{Z}_8))=sd(\mathbb{Z}_2\times Q_{16})=\frac{333}{361};$
\\

As a consequence of \textbf{Theorem 2.6}, we finish this section by providing another example of a class of groups whose subgroup commutativity degree vanishes asymptotically.
\\

\textbf{Corollary 2.8.} $\displaystyle \lim_{m \to\infty}sd(\mathbb{Z}_2 \times Q_{2^{m+1}})=0.$ 
   
\section{Cyclic subgroup commutativity degree of (generalized) dicyclic groups}

We apply similar methods as the ones we used in Section 2 to provide an explicit result for the cyclic subgroup degree of (generalized) dicyclic groups. 

Again, we start by finding the formula concerning the dicyclic group $Dic_{4n}$ and the first step is to write down its cyclic subgroups. For each divisor $r$ of $2n$, the dicyclic group has one cyclic subgroup isomorphic to $\mathbb{Z}_r$, namely $H_0^r=\langle a^{\frac{2n}{r}}\rangle$. Also, the subgroups $H_i^1=\langle a^n, a^{i-1}\gamma\rangle$, where $i=\overline{1,n}$, are cyclic since $(a^{i-1}\gamma)^2=a^n, \forall i=\overline{1,n}$. 
\\

\textbf{Theorem 3.1.}
\textit{Let $n\geqslant 1$ be a positive integer. The cyclic subgroup commutativity degree of the dicyclic group $Dic_{4n}$ is given by
$$csd(Dic_{4n})=\begin{cases} \frac{\tau (2n)(\tau (2n)+n)+n(\tau (2n)+1)}{(\tau(2n)+n)^2}, &\mbox{if } n\equiv 1(mod \ 2) \\  \frac{\tau (2n)(\tau (2n)+n)+n(\tau (2n)+2)}{(\tau(2n)+n)^2}, &\mbox{if } n\equiv 0(mod \ 2) \end{cases}.$$}
\\

\textbf{Proof.}
As we state in the Introduction, the cyclic subgroup commutativity degree of the dicyclic group $Dic_{4n}$ is
\begin{equation}\label{r2}
csd(Dic_{4n})=\frac{1}{|L_1(Dic_{4n})|}\bigg(\sum_{r|2n}|C_1(H_0^r)|+\sum\limits_{i=1}^n|C_1(H_i^1)|\bigg).
\end{equation}  
Since we specified which are the cyclic subgroups of $Dic_{4n}$, we can count them and obtain that 
$$|L_1(Dic_{4n})|=\tau(2n)+n.$$
In the proof of \textbf{Theorem 2.1.} we remarked that for each divisor $r$ of $2n$, the subgroup $H_0^r$ commutes with all subgroups of $Dic_{4n}$. It follows that 
$$\sum_{r|2n}|C_1(H_0^r)|=\tau(2n)|L_1(Dic_{4n})|=\tau(2n)(\tau(2n)+n).$$
Let $i\in \lbrace 1, 2, \ldots, n\rbrace$ and consider the cyclic subgroup $H_i^1$ of $Dic_{4n}$. We have 
$$|C_1(H_i^1)|=\bigg( \bigcup_{r|2n}\lbrace H_0^r\rbrace \bigg) \cup \bigg\lbrace H_j^1 | \ H_i^1H_j^1=H_j^1H_i^1, where \ j\in \lbrace 1, 2, \ldots, n \rbrace \bigg\rbrace.$$
Using the isomorphism $\frac{Dic_{4n}}{Z(Dic_{4n})}\cong D_{2n}$, we infer that the subgroups $H_i^1$ and $H_j^1$ of $Dic_{4n}$ commute if and only if the subgroups $K_i^1$ and $K_j^1$ of $D_{2n}$ commute. In \cite{21}, we proved that $K_j^1 \in C(K_i^1)$ iff $j=i$ and $n$ is odd or $j\in \lbrace i,(i+\frac{n}{2})(mod \ n)\rbrace$ and $n$ is even. Hence    
$$\sum\limits_{i=1}^n|C_1(H_i^1)|=\begin{cases} n(\tau(2n)+1), &\mbox{if } n\equiv 1(mod \ 2) \\  n(\tau (2n)+2), &\mbox{if } n\equiv 0(mod \ 2) \end{cases}.$$  

Making the replacements in \eqref{r2}, we obtain the desired formula and finish our proof.
\hfill\rule{1,5mm}{1,5mm}
\\

\textbf{Example 3.2.}

$\bullet \ csd(Dic_4)=csd(\mathbb{Z}_4)=1;$\\

$\bullet \ csd(Dic_8)=csd(Q_8)=1;$\\

$\bullet \ csd(Dic_{12})=\frac{43}{49};$\\

$\bullet \ csd(Dic_{16})=\frac{7}{8};$
\\

To determine an explicit formula for the cyclic subgroup commutativity degree of the generalized dicyclic group $Dic_{4n}(A)$, we need to analyse the same three cases we mentioned in \textbf{Theorem 2.6.}
\\

\textbf{Theorem 3.3.} \textit{Let $A=\mathbb{Z}_2\times \mathbb{Z}_n$ be an abelian group, where $n=2^mm'$, $m\in\mathbb{N}^*$ and $m'$ is a positive odd integer. Then, \\ \\
(i) \  if $m\geqslant 2, m'\not=1$ and $\gamma^2\in\lbrace a^{\frac{n}{2}}, b, a^{\frac{n}{2}}b\rbrace$ or if $m\geqslant 2, m'=1$ and $\gamma^2 \in \lbrace b, a^\frac{n}{2}b \rbrace$, the cyclic subgroup commutativity degree of the generalized dicyclic group $Dic_{4n}(A)$ is 
$$csd(Dic_{4n}(A))=\frac{|L_1(A)|(|L_1(A)|+n)+n(|L_1(A)|+2)}{(|L_1(A)|+n)^2}.$$\\
(ii) \ if $m\geqslant 2, m'=1$ and $\gamma^2=a^{\frac{n}{2}}$, the generalized dicyclic group $Dic_{4n}(A)$ is isomorphic to the direct product $\mathbb{Z}_2\times Q_{2^{m+1}}$ and its cyclic subgroup commutativity degree is
$$csd(\mathbb{Z}_2\times Q_{2^{m+1}})=csd(\mathbb{Z}_2\times Q_{2^{m+1}})=\frac{2^m(4m+8)+(2m+2)^2}{(2^m+2m+2)^2}.$$\\
(iii) if $m=m'=1$ and $\gamma^2\in \lbrace a^\frac{n}{2}, b, a^\frac{n}{2}b\rbrace$, the generalized dicyclic group $Dic_{4n}(A)$ is isomorphic to the abelian group $\mathbb{Z}_2^3$ and its cyclic subgroup commutativity degree is $$csd(\mathbb{Z}_2^3)=1.$$ }
\\

\textbf{Proof.}  
(i) Let $n=2^mm'$, where $m\geqslant 2$, $m'$ is a positive odd integer and $\gamma^2 \in \lbrace b, a^\frac{n}{2}b\rbrace$. The cyclic subgroups of $Dic_{4n}(A)$ are those contained in the abelian group $A$ and the subgroups $H_i^1=\langle a^{i-1}\gamma\rangle$, where $i=\overline{1,n}$. This leads us to 
$$|L_1(Dic_{4n}(A))|=|L_1(A)|+n.$$

We showed that $\frac{Dic_{4n}(A)}{\langle\gamma^2\rangle}\cong D_{2n}$ by proving the statement (i) of  \textbf{Theorem 2.6.} Through this isomorphism, in the proof of \textbf{Theorem 3.1.}, we saw that for each $i\in \lbrace 1, 2, \ldots, n\rbrace$, the subgroup $H_i^1$ has a correspondant subgroup $K_i^1=\langle x^{i-1}y\rangle$ of the dihedral group $D_{2n}$ which commutes with 2 cyclic subgroups of the same group. Also, since all cyclic subgroups of the generalized dicyclic group contained in $A$ commute with all cyclic subgroups of $Dic_{4n}(A)$, we conclude that
$$csd(Dic_{4n}(A))=\frac{|L_1(A)|(|L_1(A)|+n)+n(|L_1(A)|+2)}{(|L_1(A)|+n)^2}.$$  

The same reasoning is used if $\gamma^2=a^\frac{n}{2}$, the only change being that besides the cyclic subgroups of the generalized dicyclic group contained in $A$, we have other $n$ cyclic subgroups, namely  $\langle (ab)^{i-1}\gamma\rangle$, where $i=\overline{1,n}.$
\\

(ii) Let $m\geqslant 2, m'=1$ and $\gamma^2=a^{\frac{n}{2}}$. Denote by $n'$ the quantity $2^{m-1}$. A way to show the isomorphism $Dic_{4n}(A)\cong \mathbb{Z}_2\times Q_{2^{m+1}}$ was illustrated in \textbf{Theorem 2.6.} To find all the cyclic subgroups of $\mathbb{Z}_2 \times Q_{2^{m+1}}$, we need to view the sets $H_k, k=\overline{1,5}$ which were used to prove (ii) of the same theorem. We remark that besides the cyclic subgroups contained in $A$, the direct product $\mathbb{Z}_2 \times Q_{2^{m+1}}$ has $n'$ cyclic subgroups contained in $H_1$, namely $\langle x^{i-1}y\rangle$, and $n'$ cyclic subgroups contained in $H_3$, namely $\langle x^{i-1}yb_1\rangle,$ where $i=\overline{1,n'}$. Consequently 
$$|L_1(\mathbb{Z}_2\times Q_{2^{m+1}})|=|L_1(A)|+2n'.$$ 

To see how many pairs of cyclic subgroups contained in the sets $H_{k_1}\times H_{k_2}$ commute, we move our attention to the case $(c1)$ that we solved in the proof of (ii) of \textbf{Theorem 2.6.} We infer that we must count the solutions of the congruence $i\equiv j(mod \ \frac{n'}{2})$, where $i,j=\overline{1,n'}$. It is easy to see that for each $i$, this congruence has 2 solutions. Hence, there are $2n'$ pairs of subgroups which commute contained in each set $H_{k_1}\times H_{k_2}$, where $k_1, k_2\in \lbrace 1,3\rbrace.$ Besides them, since the cyclic subgroups of $\mathbb{Z}_2\times Q_{2^{m+1}}$ contained in $A$ commute with all cyclic subgroups of the direct product, we get other $|L_1(A)|(|L_1(A)|+4n')$ pairs. Then the cyclic subgroup commutativity degree of $\mathbb{Z}_2\times Q_{2^{m+1}}$ is
\begin{equation} \label{r5}
csd(\mathbb{Z}_2\times Q_{2^{m+1}})=\frac{|L_1(A)|(|L_1(A)|+2^{m+1})+2^{m+2}}{(|L_1(A)|+2^m)^2}.
\end{equation}

According to \textbf{Theorem 12} of \cite{19}, the number of cyclic subgroups of $A=\mathbb{Z}_2\times \mathbb{Z}_{2^m}$ is 
$$|L_1(A)|=2m+2.$$ 
Making the replacements in \eqref{r5}, we obtain
$$csd(\mathbb{Z}_2\times Q_{2^{m+1}})=\frac{2^m(4m+8)+(2m+2)^2}{(2^m+2m+2)^2}.$$ 
\\

(iii) The conclusion follows easily since $\mathbb{Z}_2^3$ is an abelian group. Our proof is complete.
\hfill\rule{1,5mm}{1,5mm}
\\

\textbf{Example 3.4.}

$\bullet \ csd(Dic_{4n}(\mathbb{Z}_2\times \mathbb{Z}_4))=csd(\mathbb{Z}_2\times Q_8)=1;$\\

$\bullet \ csd(Dic_{4n}(\mathbb{Z}_2\times \mathbb{Z}_6))=\frac{43}{49};$\\

$\bullet \ csd(Dic_{4n}(\mathbb{Z}_2\times \mathbb{Z}_8))=csd(\mathbb{Z}_2\times Q_{16})=\frac{224}{256};$
\\

Our final result states that the cyclic subgroup commutativity degree of $\mathbb{Z}_2\times Q_{2^{m+1}}$ tends to $0$ as the order of the group tends to infinity.
\\ 

\textbf{Corollary 3.5.} $\displaystyle \lim_{m \to\infty}csd(\mathbb{Z}_2 \times Q_{2^{m+1}})=0.$
\\

\section{Further research}

There is a growing interest for the probabilistic aspects of finite groups theory and the subgroup commutativity degree and the cyclic subgroup commutativity degree can still provide useful and interesting results. We indicate five open problems which may constitute the subject of some further research.
\\      

\textbf{Problem 5.1.} For a fixed $\alpha\in (0,1)$, describe the structure of finite groups $G$ satisfying $sd(G)=(\leqslant, \geqslant)\alpha$ or $csd(G)=(\leqslant, \geqslant)\alpha$. What can be said about two finite groups having the same (cyclic) subgroup commutativity degree?\\

\textbf{Problem 5.2.} Which are the connections between the (cyclic) subgroup commutativity degree of a finite group and the (cyclic) subgroup commutativity degrees of its subgroups or quotients?\\

\textbf{Problem 5.3.} Let $G$ be a finite group. Study the properties of the maps $sd:L(G)\longrightarrow [0,1], H\mapsto sd(H)$ and $csd:L(G)\longrightarrow [0,1], H\mapsto csd(H)$.\\

\textbf{Problem 5.4.} Find explicit formulas for the (cyclic) subgroup commutativity degree of the direct product $\mathbb{Z}_2\times G$, where $G$ is a group belonging to a "natural" class of groups different from the class of generalized quaternion groups. What can be said about the asymptotic behaviour of $sd(\mathbb{Z}_2\times G)$ and $csd(\mathbb{Z}_2\times G)$?\\

\textbf{Problem 5.5.} Compute explicitly the (cyclic) subgroup commutativity degree of the generalized dicyclic group $Dic_{4n}(A)$, where $A$ is an arbitrary abelian group of order $2n$.\\

\vspace*{3ex}
\small

\begin{minipage}[t]{7cm}
Marius T\u arn\u auceanu \\
Faculty of  Mathematics \\
"Al.I. Cuza" University \\
Ia\c si, Romania \\
e-mail: {\tt tarnauc@uaic.ro}
\end{minipage}
\hfill
\begin{minipage}[t]{7cm}
Mihai-Silviu Lazorec \\
Faculty of  Mathematics \\
"Al.I. Cuza" University \\
Ia\c si, Romania \\
e-mail: {\tt mihai.lazorec@student.uaic.ro}
\end{minipage}

\end{document}